\RequirePackage{fix-cm}
\documentclass[referee,envcountsect]{svjour3}

\usepackage{amsfonts}
\smartqed
\usepackage{bbding}
\usepackage{graphicx}
\usepackage{latexsym, bm, amsmath,amssymb}
\usepackage{booktabs}
\usepackage[rotateright]{rotating}
\usepackage{multirow}
\usepackage{cases}
\usepackage[misc]{ifsym}
\usepackage{algorithm,algorithmic}
\usepackage{amssymb}
\usepackage{caption}
\usepackage{mathrsfs}
\usepackage{authblk}
\usepackage{geometry}
\usepackage{graphicx}
\usepackage{ntheorem}
\allowdisplaybreaks[3]

\theorembodyfont{\upshape}
\newtheorem{lem}{Lemma}[section]
\newtheorem{thm}{Theorem}[section]

\newtheorem{cor}{Corollary}[section]
\newtheorem{alg}{Algorithm}[section]

\newtheorem{defi}{Definion}[section]
\usepackage{setspace}
\usepackage[labelfont=bf,labelsep=quad]{caption}
\def\R{{\mathbb R}}
\newcommand{\reff}[1]{(\ref{#1})}
\newtheorem{exam}{Example}

\begin{document}

\title{On Solving a Class of  Linear Semi-Infinite Programs by the Trigonometric Moment}
\author{Y. Xu  \and J. Desai \and X. Yan}
\institute{Y. Xu  (\Envelope) \at
Department of Mathematics,
Southeast University, 2 Sipailou, Nanjing, Jiangsu Province, 210096 P. R. China  \email{yi.xu1983@hotmail.com}
   \\
  J. Desai   \at Manufacturing and Industrial Engineering Cluster,
School of Mechanical and Aerospace Engineering,
Nanyang Technological University,
Singapore 639798. \email{jdesai@ntu.edu.sg}
   \\
X. Yan \at Department of Mathematics,
Taiyuan Normal University,
 Taiyuan 030012, Shanxi Province, 306019, P.R. China, \email{xihong1@e.ntu.edu.sg}
}

\date{}
\maketitle
\begin{abstract}
In this paper, we propose a new easily implementable method for solving a class of semi-infinite programs, where an approximate linear semidefinite program is constructed for the concerned semi-infinite program based on the duality theory of semi-infinite programs and the theory of the trigonometric moments. We obtain  a  $\displaystyle\frac{\ln K}{K}-$optimal solution for the semi-infinite program, where $K$ is the truncation factor. Moreover, we present some numerical examples to show the efficiency of the new method.
\end{abstract}
\keywords{Semi-Infinite Program \and Trigonometric Moment \and Semidefinite Program \and FFT}
\subclass{90C34 \and 42A70}
\section{Introduction}\label{sec0}
In this paper, we focus on designing a new method to solve semi-infinite programs (SIP). There are two main approaches that have adopted in the literature for solving semi-infinite programs, namely {\it integration} and {\it discretization} \cite{sipb1,sipb2}. In the integration-based approach, the semi-infinite constraint is replaced with the equivalent integration form. Usually, since a closed-form solution for the integral  is too difficult to be obtained for generic functions, some numerical approximation-based methods are often prescribed to compute the integral \cite{sipb3}.

An alternative set of algorithms have  been developed based on {\it discretization methods}, which include simplex-like methods and cutting plane methods \cite{sip3}. A necessary feature of such discretization methods is the computation of { interpolation points}. However, since the subprograms from solving the semi-infinite program are usually not convex programs, determining the proper interpolation points is a difficult undertaking in practice. Hence, from a computational viewpoint, integration and discretization-based approaches all have disadvantages. In addition to these two types of methods, other approaches, such as interior point methods \cite{sip4}, also suffer from similar drawbacks.

In \cite{lmb,lmp}, Lasserre discussed a moment program for solving a polynomial program.  We notice that the dual of a semi-infinite program is also a moment program. Hence, in this paper, we try to adopt a novel approach for solving a semi-infinite program  based on the theory of the {trigonometric moment}, which could alleviates all of the aforementioned drawbacks.

The remainder of this paper is organized as follows. We  present our motivation, contributions, and notations in Section \ref{sec1}, and then Section \ref{sec2} introduces the trigonometric moment problem and an associated lemma that serves as a foundation to recast this problem as a semidefinite program. Then, in Section \ref{sec3}, a detailed step-wise analysis is provided to transform the trigonometric moment program into an approximate semidefinite program. Moreover, in Section \ref{sec4},
 working under some mild assumptions, we will prove that the optimal solution of the semidefinite program got in Section \ref{sec3} is just  the approximate optimal solution of the original semi-infinite program.  Some preliminary computational results on standard test problems taken from the literature are presented in Section \ref{sec5}, along with associated insights, and finally, Section \ref{sec6} concludes this paper and presents some directions for future research.

\section{Motivation, Contributions, and Notations}\label{sec1}
Consider the following linear {\it semi-infinite program} ({\bf SIP}) and its dual program \cite{sipb1}:
\begin{subequations}\label{eq-1}
  \begin{eqnarray}
   (P)~\min & c^{\top}x &  \\
    s.t. &\displaystyle \sum_{j=1}^n a_j(t)x_j\leq a_0(t),& \forall t\in [0, 2\pi], \label{eq-1b}
  \end{eqnarray}
  \end{subequations}
%

\begin{subequations}\label{eq-2}
  \begin{eqnarray}
   (D)~ \max & -\displaystyle\int_0^{2\pi} a_0(t)\mu(dt) & \label{eq-2a}\\
    s.t. &\displaystyle\int_0^{2\pi} a_j(t)\mu(dt)=-c_j,& j=1,\cdots,n, \label{eq-2b}\\
    & \mu\in\mathbf{P}([0,2\pi]), & \label{eq-2c}
    \end{eqnarray}
  \end{subequations}
where $a_j(t)\in \mathfrak{F}([0,2\pi]), j=0,\ldots,n$, $\mathfrak{F}([0, 2\pi])$ denotes a set of functions defined on $[0, 2\pi]$, $\mathbf{P}([0,2\pi])$ is a set of  finite positive Borel measures on $[0, 2\pi].$ It is easy to see that Program ({\it D}) is essentially a {\it moment program},  and if the optimal measure $\mu$ of Program ({\it D}) is known, we can then successfully recover the optimal value of Program ({\it P}) from the solution of Program ({\it D}).

In \cite{lmb,lmp}, Lasserre discussed the following moment program:
\begin{subequations}\label{eq-3}
\begin{eqnarray}
    \min & \displaystyle\int a_0(t)\sigma(dt) &   \\
    s.t. & \sigma\in \mathbf{P}(T),
  \end{eqnarray}
  \end{subequations}
where $T=\{t : a_j(t)\geq 0, j=1,\ldots,n\}$ is a semi-algebraic set and $\mathbf{P}(T)$ is the space of finite Borel
measures on $T$.  Lasserre \cite{lmb,lmp} concluded  that the optimal value of Program \reff{eq-3} is equal to the optimal value of
\begin{subequations}\label{eq-4}
  \begin{eqnarray}
    \min & a_0(t) &\\
    s.t. &a_j(t)\geq 0,& j=1,\cdots,n,
  \end{eqnarray}
    \end{subequations}
 where $a_j(t),\  j=0,\ldots,n$ are all {polynomials}.
For solving the moment Program \reff{eq-3}, Lasserre used the {\it theory of moments} and recasted Program \reff{eq-3} into a series of successive linear semidefinite programs.

Recognizing that Programs ({\it D}) and \reff{eq-3} are both moment programs, in this research, we also apply a moment-based approach for solving Program ({\it D}), leading to the solution of Program ({\it P}). However, our problem is intrinsically very different from the ones considered by Lassere. First and foremost,  Lasserre's method is based on the theory of the {\it Hamburger moment problem}, wherein $\sigma \in \mathbf{P}(\R^n)$, while the corresponding moment $\mu$ in Program ({\it D}) only belongs to $\mathbf{P}([0,2\pi])$, which cannot be cast as a Hamburger moment problem.  Second, as $a_j(t)$ can be any generic function and not necessarily a polynomial, we cannot apply Lasserre's method directly in this case. While we could certainly use a polynomial function $b_j(t)$ to approximate $a_j(t)$, possibly by using Taylor series expansions; or the least square method; or Newton's interpolation method, such an approach inherently has several disadvantages. To note a few, Taylor series expansions require a high order of differentiability for the involved functions $a_j(t)$; the underlying integrations that exist in a least squares method are computationally very expensive; and Newton's interpolation method fails to retain and exploit the nice linear structure of Program ({\it P}).

Hence, in this paper, we adopt a novel approach for solving Program ({\it D}) based on the theory of the {trigonometric moment} that not only alleviates all of the aforementioned drawbacks, but also provides a very viable approach for solving large-scale SIP. We transform the trigonometric moment program into a linear semidefinite program. Moreover, we prove that the optimal solution of this linear semidefinite program is a  $\displaystyle\frac{\ln K}{K}-$optimal solution of the original semi-infinite program, where $K$ is the truncation factor.


In closing this section, we introduce some basic notations and related definitions that will be used in our subsequent analysis.

{\bf Notations and Symbols:}

Let $\mathbb{R}^{n\times n}$ and $\mathbb{C}^{n\times n}$ denote the set of all $n\times n$ real and complex matrices, respectively, and let $\mathbb{RS}^n$ and $\mathbb{CS}^n$ represent the set of all $n\times n$ real symmetric and complex conjugate matrices, respectively. We write  $A\succeq 0$ if $A$ is {positive semidefinite}. Let $\mathbb{CS}^n_+=\{A\in \mathbb{CS}^n: A\succeq 0\}$ and $\mathbb{RS}^n_{+}=\{A\in \mathbb{RS}^n:A\succeq 0\}$  denote the set of complex conjugate and real symmetric positive semidefinite matrices, respectively. We use $\mathcal{C}$ to denote a matrix set $\{\mathcal{C}_j\}$, $\mathcal{C}_j\in \mathbb{R}^{n\times n}$, and
denote $\mathcal{C}X$ as a vector, whose element $(\mathcal{C}X)_i=Tr(\mathcal{C}_iX)$, where $Tr(A)=\sum {A}_{j,j}$; $\mathcal{C}^Tx$ is defined as $\sum \mathcal{C}_jx_j.$

We set $i$ to represent the {\it imaginary unit}, such that $i^2=-1$; $\xi^{T}$ means the conjugate transpose of $\xi$; $\bar{\alpha}$ denotes the conjugate of  complex number $\alpha$; $||\cdot||$ is 2-norm;
 $B(x^*,\rho)$ is a hyperball of radius $\rho$ centered at $x^*$, i.e., $\{x\in \R^n:||x-x^*||\leq \rho\}$; $||\Omega||=\displaystyle\max_{x\in \Omega}||x||;$ $d(x^*,\Omega)=\min\{||x-x^*||:x\in \Omega\};$ $\mathbb{P}_{\Omega}(x^*)=\arg\min\{||x-x^*||:x\in \Omega\};$ $int(\Omega)$ is the interior set of $\Omega$; $l(a,b)$ represents the line connecting points $a$ and $b$; $S(P), F(P), $ and $v(P)$ are the optimal solution set, the feasible set, and the optimal value of Program ({\it P}), respectively.

\section{The Trigonometric Moment Problem}\label{sec2}

In a trigonometric moment problem \cite{trm1,trm2}, given a finite sequence $\{\alpha_0,\ldots,\alpha_n\}$, we are interested in determining whether there exists a positive Borel measure $\mu$, defined on the interval $[0,2\pi]$, such that
$$\alpha_k=\frac{1}{2\pi}\int_{0}^{2\pi}e^{-ikt}\mu(dt), \ k = 0,\ldots,n.$$
In other words, an affirmative answer to this moment problem would imply that $\{\alpha_0,\cdots,\alpha_n\}$ are the first $n + 1$ Fourier coefficients of some positive Borel measure $\mu$ on $[0, 2\pi]$. In fact, it has been demonstrated in \cite{trm1,trm2} that the trigonometric moment problem does indeed have a solution if the condition in the following lemma holds true.

\begin{lem}\label{lem-1}
The trigonometric moment problem has a solution, i.e., $\{\alpha_k\}_{k=0}^{n}$ is a valid sequence of Fourier coefficients of some positive Borel measure $\mu$ defined on $[0, 2\pi]$, if and only if there exists a Toeplitz matrix $A\in \mathbb{CS}^{n+1}$ that is positive semidefinite, that is

$$A=\left(
   \begin{array}{cccc}
     \alpha_0 & \alpha_1 & \cdots & \alpha_n \\
     \bar{\alpha}_1 & \alpha_0 & \cdots & \alpha_{n-1} \\
     \vdots & \vdots & \ddots & \vdots \\
     \bar{\alpha}_n & \bar{\alpha}_{n-1} & \cdots & \alpha_0 \\
   \end{array}
 \right)\succeq 0.$$

 \end{lem}

\begin{remark}\label{Re1}
(a) Note that the existence of a solution for the trigonometric moment problem is similar to the one of the Hamburger moment problem, which also requires matrix $A$ to be positive semidefinite. However, the matrix $A$ in the Hamburger moment problem is a Hankel matrix, and not a Toeplitz matrix, as needed for the problem under consideration.


 (b) Matrix $A$ in Lemma \ref{lem-1} is not only a Toeplitz matrix, but it is also a {conjugate} matrix. Hence, when $A\succeq 0$, $\alpha_0$ must be a real and nonnegative number.
\end{remark}

\section{Transforming a Linear Semi-Infinite Program into an Approximate  Semidefinite Program}\label{sec3}
%

%
%

The foregoing Dirichlet conditions of $f(t)$ on $[0, 2\pi]$:
\begin{itemize}
  \item[(1)] $f(t)$ is absolutely integrable over $[0, 2\pi]$;
  \item[(2)] $f(t)$ has a finite number of extrema in $[0,2\pi]$;
  \item[(3)] $f(t)$ has a finite number of discontinuities in $[0, 2\pi]$.
\end{itemize}
 are well-known sufficiency conditions for the existence of a Fourier transformation of $f(t)$. Working under Dirichlet conditions of $a_j(t), j=0,\cdots,n$ on $[0, 2\pi]$, $a_j(t)$ can be represented via its Fourier transform as $$\hat{a}_j(t)=\frac{1}{2\pi}\sum_{k=-\infty}^{+\infty} a_{j,k}e^{-ikt},$$
where
\begin{eqnarray}\label{F1}
a_{j,k}=\displaystyle\frac{1}{2\pi}\int_0^{2\pi}a_j(t)e^{ikt}dt.
\end{eqnarray}

\noindent
Substituting $\hat{a}_j(t)$ in lieu of $a_j(t)$ into Program ({\it D}), the resulting \textit{Fourier Dual}  Program ({\it FD}) can be written as:
 \begin{subequations}\label{eq-3.1}
  \begin{eqnarray}
   (FD) ~\max & -\displaystyle\frac{1}{2\pi}\int_0^{2\pi} \sum_{k=-\infty}^{+\infty}a_{0,k}e^{-ikt}\mu(dt) &\\
    s.t. &\displaystyle\frac{1}{2\pi}\int_0^{2\pi} \sum_{k=-\infty}^{+\infty}a_{j,k}e^{-ikt}\mu(dt)=-c_j,& j=1,\ldots,n,\\
     & \mu\in\mathbf{P}([0,2\pi]). &
  \end{eqnarray}
  \end{subequations}
As Program ($FD$) cannot be directly solved due to the presence of the infinite summation terms present in this program, we consider a {truncated} version of this program, denoted as Program ($FD_K$), given by:
 \begin{subequations}\label{eq-3.5}
  \begin{eqnarray}
  (FD_K)~  \max & -\displaystyle\frac{1}{2\pi}\int_0^{2\pi} \sum_{k=-K}^{K}a_{0,k}e^{-ikt}\mu(dt) &\\
    s.t. &\displaystyle\frac{1}{2\pi}\int_0^{2\pi} \sum_{k=-K}^{K}a_{j,k}e^{-ikt}\mu(dt)=-c_j, &  j=1,\ldots,n,\\
    & \mu\in\mathbf{P}([0,2\pi]), &
  \end{eqnarray}
    \end{subequations}
  where $K$, a positive integer, is the truncation factor.

   Interchanging the summation and integral terms, this program can in turn be represented as:
  \begin{subequations}\label{eq-4}
    \begin{eqnarray}
  (FD_K)~  \max & -\displaystyle\sum_{k=-K}^{K}a_{0,k}\frac{1}{2\pi}\int_0^{2\pi} e^{-ikt}\mu(dt) &\\
    s.t. &\displaystyle\sum_{k=-K}^{K}a_{j,k} \frac{1}{2\pi}\int_0^{2\pi} e^{-ikt}\mu(dt)=-c_j,&  j=1,\cdots,n,\label{conX}\\
     & \mu \in\mathbf{P}([0,2\pi]). &
  \end{eqnarray}
  \end{subequations}

Now, let $y_k=\displaystyle\frac{1}{2\pi}\displaystyle\int_0^{2\pi} e^{-ikt}\mu(dt),\ k=-K,\ldots,K$. It is easy to see that $y_0$ is a real number, in general,  $y_k$ satisfies $y_{-k}=\bar{y}_{k}$.

 From Lemma \ref{lem-1}, it follows that the matrix composed of  $y_k$ satisfies
$$\left(
   \begin{array}{cccc}
     y_0 & y_1 & \cdots & y_K \\
     \bar{y}_1 & y_0 & \cdots & y_{K-1} \\
     \vdots & \vdots & \ddots & \vdots \\
     \bar{y}_K & \bar{y}_{K-1} & \cdots & y_0 \\
   \end{array}
 \right)=\left(
   \begin{array}{cccc}
     y_0 & y_1 & \cdots & y_K \\
     y_{-1} & y_0 & \cdots & y_{K-1} \\
     \vdots & \vdots & \ddots & \vdots \\
     {y}_{-K} & {y}_{-K+1} & \cdots & y_0 \\
   \end{array}
 \right)\succeq 0.$$

Conversely, once again by Lemma \ref{lem-1}, if $y_k$ satisfies the above positive semidefiniteness condition, then, there must be a positive Borel measure $\mu$, defined on the interval $[0,2\pi]$, such that  $$y_k=\frac{1}{2\pi}\displaystyle\int_0^{2\pi} e^{-ikt}\mu(dt),\ k=-K,\ldots,K.$$

Incorporating the results of Lemma \ref{lem-1} into Program $(FD_K)$, we get the following {\it Trigonometric Fourier Dual} Program $(TFD_{K})$ variant of Program ({\it P}):

   \begin{subequations}\label{eq-5}
    \begin{eqnarray}
  (TFD_K)~  \max & -\displaystyle\sum_{k=-K}^{K} a_{0,k} y_k  \\
    s.t. &\displaystyle\sum_{k=-K}^{K}a_{j,k} y_k=-c_j, j=1,\cdots,n,\label{conX1}\\
       & \left(
   \begin{array}{cccc}
     y_0 & y_1 & \cdots & y_K \\
     y_{-1} & y_0 & \cdots & y_{K-1} \\
     \vdots & \vdots & \ddots & \vdots \\
     {y}_{-K} & {y}_{-K+1} & \cdots & y_0 \\
   \end{array}
 \right)\succeq 0.
  \end{eqnarray}
\end{subequations}

 It is clear that Program $(TFD_K)$ is not a conventional optimization problem as both $a_{j,k}$ and $y_k$ in Program $(TFD_K)$ may be complex numbers, for all $k\neq 0$. Hence, for solving Program $(TFD_K)$, we adopt a novel approach by separating each of the program parameters and variables into their respective real and imaginary components, and then using the properties of the involved matrices to derive an equivalent program defined purely in terms of real-valued variables. Towards this end, set
\begin{subequations}
\begin{eqnarray}
y_k&=&w_k+iv_k,		\label{eq-9a} \\
a_{j,k}&=&r_{j,k}+i s_{j,k}, 	\label{eq-9b}
\end{eqnarray}
\end{subequations}
where $w_k, v_k, r_{j,k}$, and $ s_{j,k}$ are all real numbers.

Then,
\begin{subequations}
\begin{eqnarray}\label{eq-10}
&&w_k =w_{-k},  v_k=-v_{-k}, v_0=0 \\
 &&r_{j,k}=r_{j,-k},  s_{j,k}=-s_{j,-k},\label{eq-10b}\\
 &&   s_{j,0}=0, j=0,\ldots,n, \label{eq-10c}
\end{eqnarray}
\end{subequations}
  where equations \eqref{eq-10b} and \eqref{eq-10c} are  a byproduct of the definition of the Fourier transformation of $a_j(t)$.

Substituting these resulting terms into Program $(TFD_K)$, we get

$$\displaystyle\sum_{k=-K}^{K}a_{j,k}y_k=\displaystyle\sum_{k=-K}^{K}r_{j,k}w_k-s_{j,k}v_k, \  j=0,\ldots,n,$$
and $$\left(
   \begin{array}{cccc}
     y_0 & y_1 & \cdots & y_K \\
     y_{-1} & y_0 & \cdots & y_{K-1} \\
     \vdots & \vdots & \ddots & \vdots \\
     {y}_{-K} & {y}_{-K+1} & \cdots & y_0 \\
   \end{array}
 \right)=\left(
   \begin{array}{cccc}
     w_0 & w_1 & \cdots & w_K \\
     w_{1} & w_0 & \cdots & w_{K-1} \\
     \vdots & \vdots & \ddots & \vdots \\
     {w}_{K} & {w}_{K-1} & \cdots & w_0 \\
   \end{array}
 \right) + i\left(
   \begin{array}{cccc}
     0 & v_1 & \cdots & v_K \\
     -v_{1} & 0 & \cdots & v_{K-1} \\
     \vdots & \vdots & \ddots & \vdots \\
     -{v}_{K} & -{v}_{K-1} & \cdots & 0 \\
   \end{array}
 \right)\succeq 0. $$

We are now ready to present Lemma \ref{lem-2}, which accords us with a mechanism to transform $Y\in \mathbb{CS}^{K}, Y\succeq 0$ into a real-valued positive semidefinite matrix.

 \begin{lem}\label{lem-2}
$$\left(
   \begin{array}{cccc}
     y_0 & y_1 & \cdots & y_K \\
     y_{-1} & y_0 & \cdots & y_{K-1} \\
     \vdots & \vdots & \ddots & \vdots \\
     {y}_{-K} & {y}_{1-K} & \cdots & y_0 \\
   \end{array}
 \right)\succeq 0
\Leftrightarrow\left(
   \begin{array}{cccccccc}
     w_0 & w_1 & \cdots & w_K&0 & -v_1 & \cdots & -v_K \\
     w_{1} & w_0 & \cdots & w_{K-1}& v_{1} & 0 & \cdots & -v_{K-1}\\
     \vdots & \vdots & \ddots & \vdots& \vdots & \vdots & \ddots & \vdots\\
     {w}_{K} & {w}_{K-1} & \cdots & w_0&{v}_{K} & {v}_{K-1} & \cdots & 0 \\
     0 & v_1 & \cdots & v_K &w_0 & w_1 & \cdots & w_K\\
     -v_{1} & 0 & \cdots & v_{K-1}&w_{1} & w_0 & \cdots & w_{K-1}\\
      \vdots & \vdots & \ddots & \vdots& \vdots & \vdots & \ddots & \vdots\\
      -{v}_{K} & -{v}_{K-1} & \cdots & 0&{w}_{K} & {w}_{K-1} & \cdots & w_0\\
   \end{array}
 \right)\succeq 0.$$
\end{lem}
{\noindent{\rm{ \bf Proof~}} } Denote $$Y=\left(
   \begin{array}{cccc}
     y_0 & y_1 & \cdots & y_K \\
     y_{-1} & y_0 & \cdots & y_{K-1} \\
     \vdots & \vdots & \ddots & \vdots \\
     {y}_{-K} & {y}_{-K+1} & \cdots & y_0 \\
   \end{array}
 \right), W=\left(
   \begin{array}{cccc}
     w_0 & w_1 & \cdots & w_K \\
     w_{1} & w_0 & \cdots & w_{K-1} \\
     \vdots & \vdots & \ddots & \vdots \\
     {w}_{K} & {w}_{K-1} & \cdots & w_0 \\
   \end{array}
 \right), V=\left(
   \begin{array}{cccc}
     0 & v_1 & \cdots & v_K \\
     -v_{1} & 0 & \cdots & v_{K-1} \\
     \vdots & \vdots & \ddots & \vdots \\
     -{v}_{K} & -{v}_{K-1} & \cdots & 0 \\
   \end{array}
 \right),$$
where, clearly $W = W^{\top}, V =-V^{\top}$, with $W, V\in \R^{(K+1)\times (K+1)}$.

 As $Y\succeq 0$, this implies that $(p+iq)^{T}Y(p+iq)\geq 0, \forall p, q\in \mathbb{R}^{(K+1)\times 1}$. Expanding this product, we get
  \begin{eqnarray}
   (p+iq)^{T}Y(p+iq) &=&p^{\top}Wp-p^{\top}Vq+q^{\top}Wq+q^{\top}Vp+i(p^{\top}Wq+p^{\top}Vp-q^{\top}Wp+q^{\top}Vq).\nonumber
 \end{eqnarray}
 Utilizing the property that $W$ is a symmetric matrix and $V$ is a skew-symmetric matrix, which implies $p^{\top}Wq=q^{\top}Wp$ and $p^{\top}Vp=q^{\top}Vq=0,$ we get that $Y\succeq 0$ if and only if $$p^{\top}Wp-p^{\top}Vq+q^{\top}Wq+q^{\top}Vp\geq 0, \forall p, q\in \mathbb{R}^{(K+1) \times 1},$$ which in turn can be compactly expressed as a real, symmetric  matrix
$$\left(
    \begin{array}{cc}
   W & V^{\top}\\
   V &W\\
   \end{array}
 \right)\succeq 0,$$ and this completes the proof. $\hfill \square$

Using the variable transformations \reff{eq-9a} and \reff{eq-9b}, in concert with Lemma \ref{lem-1}, the following equivalent form of Program ($TFD_K$), now defined in terms of only real-valued variables, can be obtained as follows:
\begin{subequations}\label{eq-6}
    \begin{eqnarray}
  (TFD_K)~  \max & -\displaystyle\sum_{k=-K}^{K}r_{0,k}w_k-s_{0,k}v_k & \\
    s.t. &\displaystyle\sum_{k=-K}^{K}r_{j,k}w_k-s_{j,k}v_k=-c_j,\quad  j=1,\ldots,n, &\label{eq-11b}\\
    &v_{k}=-v_{-k}, w_k=w_{-k},\quad  k=0,\ldots,K, &\\
       & \left(
   \begin{array}{cccccccc}
     w_0 & w_1 & \cdots & w_K&0 & -v_1 & \cdots & -v_K \\
     w_{1} & w_0 & \cdots & w_{K-1}& v_{1} & 0 & \cdots & -v_{K-1}\\
     \vdots & \vdots & \ddots & \vdots& \vdots & \vdots & \ddots & \vdots\\
     {w}_{K} & {w}_{K-1} & \cdots & w_0&{v}_{K} & {v}_{K-1} & \cdots & 0 \\
     0 & v_1 & \cdots & v_K &w_0 & w_1 & \cdots & w_K\\
     -v_{1} & 0 & \cdots & v_{K-1}&w_{1} & w_0 & \cdots & w_{K-1}\\
      \vdots & \vdots & \ddots & \vdots& \vdots & \vdots & \ddots & \vdots\\
      -{v}_{K} & -{v}_{K-1} & \cdots & 0&{w}_{K} & {w}_{K-1} & \cdots & w_0\\
   \end{array}
 \right)\succeq 0.
   \end{eqnarray}
   \end{subequations}

A straightforward extension of the above discussion results in the following theorem, it is obliviously that $v(TFD_K)=v(FD_K).$

Let $$\bar{\mathcal{C}}^K_1=\{\bar{\mathcal{C}}_{1,j}^K\}_{j=0}^K, \bar{\mathcal{C}}^K_{1,j}=\left(
                                                                                          \begin{array}{cc}
                                                                                            \mathcal{C}^K_{1,j} & 0 \\
                                                                                            0 & \mathcal{C}^K_{1,j} \\
                                                                                          \end{array}
                                                                                        \right),
 (\mathcal{C}^K_{1,j})_{k,l}=\left\{
                                                                                        \begin{array}{cc}
                                                                                          1, & |k-l|=j \\
                                                                                          0, & else \\
                                                                                        \end{array}
                                                                                      \right.
,\mathcal{C}^K_{1,j}\in \mathbb{R}^{(K+1)\times(K+1)},$$  $$\bar{\mathcal{C}}^K_2=\{\bar{\mathcal{C}}^K_{2,j}\}_{j=1}^K,  \bar{\mathcal{C}}^K_{2,j}=\left(
                                                                                                                                               \begin{array}{cc}
                                                                                                                                                0& -\mathcal{C}_{2,j}^K \\
                                                                                                                                                 \mathcal{C}_{2,j}^K & 0 \\
                                                                                                                                               \end{array}
                                                                                                                                             \right)
,  (\mathcal{C}^K_{2,j})_{k,l}=\left\{
                                                                                        \begin{array}{cc}
                                                                                          1, & l-k=j \\
                                                                                          -1,& k-j=j \\
                                                                                          0, & else \\
                                                                                        \end{array}
                                                                                      \right.
, \mathcal{C}^K_{2,j}\in \mathbb{R}^{(K+1)\times(K+1)},$$ $$ p^K_j= \left\{
                                                \begin{array}{cc}
                                                  r_{0,0} & j=0 \\
                                                  r_{0,j}+r_{0,-j} & else \\
                                                \end{array}
                                              \right., p^K\in \mathbb{R}^{(K+1)\times 1},
 q^K_j=s_{0,j}-s_{0,-j}, q^K\in \mathbb{R}^{K\times 1},$$ $$E^K_{j,k}=\left\{
                                                                                                       \begin{array}{cc}
                                                                                                         r_{j, 0}, & k=0 \\
                                                                                                         r_{j, k}+r_{j,-k}, & else \\
                                                                                                       \end{array}
                                                                                                     \right., E^K\in \mathbb{R}^{n\times (K+1)},
 F^K_{j,k}=s_{j,k}-s_{j,-k}, F^K\in \mathbb{R}^{n\times K},$$ $$w=(w_0,w_1,\cdots,w_K)^T, v=(v_1,\cdots,v_K)^T,$$ $(TFD_K)$ could be rewritten as

\begin{subequations}\label{eq-7}
\begin{eqnarray}
  (TFD_K) ~-\max & (p^K)^Tw-(q^K)^Tv &\nonumber\\
        s.t. & E^Kw-F^Kv=-c, &\nonumber\\
        & (\bar{\mathcal{C}}^K_1)^Tw+(\bar{\mathcal{C}}_2^K)^Tv\succeq 0.&\nonumber
\end{eqnarray}
   \end{subequations}
Furthermore, its dual program $(TFP_K)$ is the following semidefinite program:
\begin{subequations}\label{eq-8}
\begin{eqnarray}
  (TFP_K) ~\min & c^Tx &\nonumber\\
        s.t. & \bar{\mathcal{C}}^K_1Y+(E^K)^Tx+p^K=0, &\nonumber\\
        & \bar{\mathcal{C}}^K_2Y-(F^K)^Tx-q^K=0,&\nonumber\\
        & Y\succeq 0 .\nonumber
\end{eqnarray}
   \end{subequations}

%
%
%

\section{Error Analysis and Algorithm}\label{sec4}

Firstly, we introduce the definition of $\varepsilon-$optimal solution here.
\begin{defi}
We call a point $x^*$ is a $\varepsilon-$optimal solution of $(P)$, if there is a constant $L>0$ such that
\begin{itemize}
  \item[(a)] $|c^Tx^*-v(P)|\leq L\varepsilon;$
  \item[(b)] $d(x^*,F(P))\leq L\varepsilon.$
\end{itemize}
\end{defi}

\begin{lem}\label{lem-0}
If $\exists~ \tilde{x}\in \R^n$ and $Q<0$ such that  $$\displaystyle\sum_{j=1}^na_j(t)\tilde{x}_j\leq Q, \forall t\in[0,2\pi],$$
then, there exists $\sigma>0$, only depended on $\tilde{x}$ and $Q$, such that $$d(x,F(P))\leq \sigma \max\{ \max_{t\in[0,2\pi]}\sum_{j=1}^na_j(t)x_j-a_0(t),0\}$$
\end{lem}
\proof If $x\in F(P)$, it is true, so we only need to consider $x\not\in F(P)$.
Let $y=\mathbb{P}_{F(P)}(x)$, there must be $\lambda_j\geq 0, t_j\in [0,2\pi],j=1,\cdots,l$ such that $0\neq x-y=A^T\lambda$, $Ay=b$ where $A\in \mathbb{R}^{l\times n }, A_{k,j}=a_j(t_k)$,  $b=(a_0(t_1),a_0(t_2),\cdots,a_0(t_l))^T$, so $d(x,F(P))=||x-y||$.


Then, $\max\{\displaystyle\max_{t\in [0,2\pi]}\sum_{j=1}^na_j(t)x_j-a_0(t),0\}=\max_{j=1,\cdots,l}\{\max\{(AA^T\lambda)_j,0\}\}.$
Since $$\sqrt{n}\max_{j=1,\cdots,l}\{\max\{(AA^T\lambda)_j,0\}\}||\lambda||\geq \sum_{j=1,\cdots,l}\lambda_j\max\{(AA^T\lambda)_j,0\}\geq \sum_{j=1,\cdots,l}\lambda_j(AA^T\lambda)_j=||A^T\lambda||^2,$$
then,
$$\frac{\sqrt{n}\displaystyle\max_{j=1,\cdots,l}\{\max\{(AA^T\lambda)_j,0\}\}}{||A^T\lambda||}\geq \frac{||A^T\lambda||}{||\lambda||}.$$

Since  $$\displaystyle\sum_{j=1}^na_j(t)\tilde{x}_j\leq Q, \forall t\in[0,2\pi],$$ we have $A(-\tilde{x})>0$,
so $$||\tilde{x}||||A^T\lambda||\geq -\sum_{j=1}^n\tilde{x}_j(A^T\lambda)_j\geq -\max_{j=1,\cdots,n}\{(A\tilde{x})_j\}||\lambda||\geq -Q||\lambda||,$$
then, \begin{eqnarray}
             d(x,F(P))=&&||x-y||=||A^T\lambda||\leq \frac{||\lambda||}{||A^T\lambda||}\displaystyle\sqrt{n}\max_{j=1,\cdots,l}\{\max\{(AA^T\lambda)_j,0\}\}\nonumber \\
               \leq&&\frac{\sqrt{n}||\tilde{x}||}{-Q}\max_{j=1,\cdots,l}\{\max\{(AA^T\lambda)_j,0\} \nonumber\\
               =&&\sigma\max\{\displaystyle\max_{t\in [0,2\pi]}\sum_{j=1}^na_j(t)x_j-a_0(t),0\},\nonumber
             \end{eqnarray}
             where $\sigma=\displaystyle\frac{\sqrt{n}||\tilde{x}||}{-Q}.$
\qed

Now, we are ready to present our main result.

\begin{thm}\label{thm-0}
Supposing the following conditions hold:
\begin{itemize}
  \item[(a)] $a_j(t), j=0,\cdots,n$ are Lipschitz continuous functions on $[0, 2\pi]$;
    \item[(b)] $a_j(0)=a_j(2\pi),j=0,\cdots,n$;
  \item[(c)] The slater condition holds for Program $(P)$:
there exists a point $\hat{x}$ satisfied $$\displaystyle\sum_{j=1}^na_j(t)\hat{x}_j<a_0(t), \forall t\in[0,2\pi];$$
\item[(d)]  $\exists~ \tilde{x}\in \R^n$ and $Q<0$ such that  $$\displaystyle\sum_{j=1}^na_j(t)\tilde{x}_j\leq Q, \forall t\in[0,2\pi];$$
\item[(e)] $c\neq 0,$ $S(P)\neq \emptyset$, $||S(P)||\leq R,$
\end{itemize}
 then,
\begin{itemize}
  \item[(a)] the slater condition holds for Program $(FP_K)$,  when $K$ is large enough, where
 \begin{subequations}\label{eq-3.5}
  \begin{eqnarray}
  (FP_K)~  \min & c^{\top}x &\nonumber\\
    s.t. &\displaystyle\sum_{j=1}^n \hat{a}^K_j(t)x_j\leq \hat{a}^K_0(t), & \forall t\in [0, 2\pi] \nonumber
  \end{eqnarray}
  \end{subequations}
is the dual program of Program $(FD_K)$;
  \item[(b)] if $\tilde{x}^K\in S(AFP_K)$, $\tilde{x}^K$ is a $\displaystyle \frac{\ln K}{K}-$optimal solution of Program $(P)$, where   \begin{subequations}\label{eq-3.6}
  \begin{eqnarray}
  (AFP_K)~  \min & c^{\top}x &\nonumber\\
    s.t. &\displaystyle\sum_{j=1}^n \hat{a}^K_j(t)x_j\leq \hat{a}^K_0(t)+2(nR+1)A\frac{\ln K}{K}, & \forall t\in [0, 2\pi], \nonumber
  \end{eqnarray}
  \end{subequations}
 when $K$ is large enough;
  \item[(c)] if $(\tilde{x}^*)^K\in S(FP_K)$, $(\tilde{x}^*)^K$ is a $\displaystyle \frac{\ln K}{K}-$optimal solution of Program $(P)$, when $K$ is large enough;
  \item[(d)]    $(x^{*})^K$ is a $\displaystyle\frac{\ln K}{K}-$optimal solution of Program $(P)$,   when $K$ is large enough, where $((x^{*})^K,(Y^*)^K)\in S(TFP_K)$ with some $(Y^*)^K\succeq 0.$
\end{itemize}

\end{thm}
\proof Since $a_j(t)$ is a Lipschitz function on $[0, 2\pi]$ and $a_j(0)=a_j(2\pi)$, we always could extend $a_j(t)$ to $\breve{a}_j(t)$ defined on $\R$ such that $\breve{a}_j(t)$ is a Lipschitz function on $\R$. From  Corollary I in page 22 of \cite{book2}, it follows that there exists a $A>0$ such that $$|a_j(t)-\hat{a}^K_j(t)|\leq \displaystyle A\frac{\ln K}{K}, \forall t\in[0,2\pi],$$
where $\hat{a}^K_j(t)$ is the  Fourier transformation of $a_j(t).$


(a) Let $K$ be large enough such that $$\displaystyle\min_{t\in[0,2\pi]}(a_0(t)-\sum_{j=1}^na_j(t)\hat{x}_j) -(n||\hat{x}||+1)A\frac{\ln K}{K}>0,$$
and $$\tilde{\rho}= 0.9\displaystyle\frac{\displaystyle\min_{t\in[0,2\pi]}(a_0(t)-\sum_{j=1}^na_j(t)\hat{x}_j) -(n||\hat{x}||+1)A\frac{\ln K}{K}}{\displaystyle\max_{t\in[0,2\pi]}\displaystyle\sum_{j=1}^n|a_j(t)|+nA\frac{\ln K}{K}}>0,$$
then, $\forall \check{x}\in B(\hat{x},\tilde{\rho})$,
\begin{eqnarray*}
   &&\max_{t\in[0,2\pi]}\displaystyle\sum_{j=1}^n\hat{a}^K_j(t)\check{x}_j-a^K_0(t)  \nonumber\\
   =&&\max_{t\in[0,2\pi]}\displaystyle\sum_{j=1}^n\hat{a}^K_j(t)\hat{x}_j-a^K_0(t)+\displaystyle\sum_{j=1}^n\hat{a}^K_j(t)(\check{x}_j-\hat{x}_j)\nonumber\\
   =&&\max_{t\in[0,2\pi]}\displaystyle\sum_{j=1}^n(\hat{a}^K_j(t)-a_j(t))\hat{x}_j-(a^K_0(t)-a_0(t))+\displaystyle\sum_{j=1}^n(\hat{a}^K_j(t)-a_j(t))(\check{x}_j-\hat{x}_j) \nonumber\\
   &&+\displaystyle\sum_{j=1}^na_j(t)\hat{x}_j-a_0(t)+\displaystyle\sum_{j=1}^na_j(t)(\check{x}_j-\hat{x}_j)   \nonumber\\
   \leq&&(n||\hat{x}||+1)A\frac{\ln K}{K}+\tilde{\rho}nA\frac{\ln K}{K}+\tilde{\rho}\displaystyle\max_{t\in[0,2\pi]}\sum_{j=1}^n|a_j(t)|+\max_{t\in[0,2\pi]}\displaystyle\sum_{j=1}^na_j(t)\hat{x}_j-a_0(t)  \nonumber\\
   <&&0,
\end{eqnarray*} so the slater condition holds for Program $(FP_K)$, when $K$ is large enough.

(b) It is obviously that the slater condition also holds for Program $(AFP_K)$. Furthermore, since $S(P)\neq \emptyset, ||S(P)||\leq R$, then, $\forall \bar{x}\in S(P)$,
\begin{eqnarray*}
   &&\max_{t\in[0,2\pi]}\displaystyle\sum_{j=1}^n\hat{a}^K_j(t)\bar{x}_j-\hat{a}^K_0(t)  \nonumber\\
   =&&\max_{t\in[0,2\pi]}\displaystyle\sum_{j=1}^n(\hat{a}^K_j(t)-a_j(t))\bar{x}_j-(\hat{a}^K_0(t)-a_0(t))+\displaystyle\sum_{j=1}^na_j(t)\bar{x}_j-a_0(t)\nonumber\\
   <&&2(nR+1)A\frac{\ln K}{K},  \nonumber
\end{eqnarray*}
which means $\bar{x}\in int(F(AFP_K))$ and $v(P)> v(AFP_K).$

If $||S(AFP_K)||>2R$, there exist $\{K_i\}$ such that  $||\tilde{x}^{K_i}||>2R, \tilde{x}^{K_i}\in S(AFP_{K_i}).$ Taken $\tilde{\tilde{x}}^{K_i}\in l(\bar{x},\tilde{x}^{K_i})$ satisfied $||\tilde{\tilde{x}}^{K_i}||=2R,$ where $\bar{x}\in S(P)$, we have $v(AFP_{K_i})\leq c^T\tilde{\tilde{x}}^{K_i}<v(P)$, $\tilde{\tilde{x}}^{K_i}\in F(AFP_{K_i})$.  By Lemma \ref{lem-0}, $$d(\tilde{\tilde{x}}^{K_i},F(P))\leq \sigma \max\{\max_{t\in [0,2\pi]} \sum_{j=1}^na_j(t)\tilde{\tilde{x}}^{K_i}_j-a_0(t) ,0\}.$$

Noticing that \begin{eqnarray}
              &&\max_{t\in [0,2\pi]} \sum_{j=1}^na_j(t)\tilde{\tilde{x}}^{K_i}_j-a_0(t) \nonumber\\
                 =&&  \max_{t\in [0,2\pi]} \sum_{j=1}^n(a_j(t)-\hat{a}_j^{K_i}(t))\tilde{\tilde{x}}^{K_i}_j-(a_0(t)-\hat{a}_0^{K_i}(t))+  \sum_{j=1}^n\hat{a}_j^{K_i}(t)\tilde{\tilde{x}}^{K_i}_j-\hat{a}_0^{K_i}(t)\nonumber\\
                 \leq &&(2nR+1)A\frac{\ln K_i}{K_i}+2(nR+1)A\frac{\ln K_i}{K_i} \nonumber\\
                 =&&(4nR+3)A\frac{\ln K_i}{K_i}. \nonumber
              \end{eqnarray}
       Then, $$d(\tilde{\tilde{x}}^{K_i},F(P))\leq \sigma(4nR+3)A\frac{\ln K_i}{K_i}.$$
       When $K_i\rightarrow \infty$, $d(\tilde{\tilde{x}}^{K_i},F(P)) \rightarrow 0$. For any accumulation point $\tilde{\tilde{x}}$ of $\tilde{\tilde{x}}^{K_i}$, we have  $c^T\tilde{\tilde{x}}=v(P),$ $||\tilde{\tilde{x}}||=2R$ , which is against to  $||S(P)||\leq R$. So $||S(AFP_K)||\leq 2R$, when $K$ is large enough.

By the same way, we could prove that $$d(\tilde{x}^{K},F(P))\leq \sigma(4nR+3)A\frac{\ln K}{K},$$
where $\tilde{x}^K\in S(AFP_K).$
Then, $$|c^T\tilde{x}^K-c^T\mathbb{P}_{F(P)}{\tilde{x}^K}|\leq ||c||\sigma(4nR+3)A\frac{\ln K}{K},$$
which means $$-L\frac{\ln K}{K}\leq v(AFP_K)-v(P)\leq 0,$$
where $L=||c||\sigma(4nR+3)A.$

(c) We will prove there exists two constants $L>0$ and $U>0$ such that $$-L{\frac{\ln K}{K}}\leq v(FP_K)-v(P)\leq U{\frac{\ln K}{K}}.$$

Since the slater condition holds for Program $(AFP_K)$,
\begin{eqnarray}
  0\geq &&v(AFP_K)-v(FD_K)=v(AFD_K)-v(FD_K)\nonumber\\
   \geq && \begin{array}{rlcl}
             -&\max & \displaystyle2(nR+1)A\frac{\ln K}{K}\displaystyle\int_0^{2\pi} \mu(dt) &  \\
             &s.t.  & \displaystyle\int_0^{2\pi}\hat{a}^K_j(t)\mu(dt)=-c_j, &  j=1,\cdots,n,\\
              & & \mu\in\mathbf{P}([0,2\pi]), & \nonumber
           \end{array}
      \nonumber
\end{eqnarray}
where \begin{eqnarray}
  (AFD_K)~  \max & -\displaystyle\int_0^{2\pi}\hat{a}_0^{K}(t)+2(nR+1)A\frac{\ln K}{K} \mu(dt) &\nonumber \\
    s.t. &\displaystyle\int_0^{2\pi}\hat{a}_j^{K}(t) \mu(dt)=-c_j, \quad j=1,\ldots,n,&\nonumber\\
     & \mu\in\mathbf{P}([0,2\pi]) &\nonumber
  \end{eqnarray}
is the dual program of Program $(AFP_K).$

We consider  Program $(TD)$ and its dual Program $(TP)$:
\begin{eqnarray}
 (TD)~ \max & -\displaystyle\frac{Q}{2}\displaystyle\int_0^{2\pi} \mu(dt) & \nonumber\\
  s.t. &\displaystyle\int_0^{2\pi}\hat{a}^K_j(t) \mu(dt)=-c_j,& j=1,\cdots,n, \nonumber\\
  & \mu\in\mathbf{P}([0,2\pi]), & \nonumber
  \end{eqnarray}
  \begin{eqnarray}
 (TP)~ \min & c^Tx                                &\nonumber\\
  s.t. &\displaystyle\sum_{j=1}^n\hat{a}^K_j(t)x_j\leq \frac{Q}{2}   ,& \forall t\in [0,2\pi].\nonumber
\end{eqnarray}

It is easy to see that $v(TP)\geq v(TD)\geq 0$, since $F(TD)=F(AFD_K)$ is not empty.
Let $K$ be large enough such that  $\displaystyle A\frac{\ln K}{K}\leq \displaystyle\frac{-Q}{2n||\tilde{x}||}$,
 then,
\begin{eqnarray}
  \displaystyle\sum_{j=1}^n\hat{a}^K_j(t)\tilde{x}_j=&& \displaystyle\sum_{j=1}^n(\hat{a}^K_j(t)-a_j(t))\tilde{x}_j+\displaystyle\sum_{j=1}^na_j(t)\tilde{x}_j\leq n||\tilde{x}||A\frac{\ln K}{K}+Q\leq  \frac{Q}{2},\nonumber
\end{eqnarray}
 which means that $\tilde{x}$ is a feasible point of Program $(TP)$, then,
$$0\leq v(TD)\leq v(TP)\leq c^{\top}\tilde{x}.$$ Let $U=\displaystyle\frac{-4(nR+1)c^{\top}\tilde{x}}{Q}A,$ we have that
$$-U\frac{\ln K}{K} \leq v(AFD_K)-v(FD_K)= v(AFP_K)-v(FD_K)\leq 0,$$
since the slater condition holds for Program $(AFP_K)$ , when $K$ is large enough. So when $K$ is large enough, we have
$$-L{\frac{\ln K}{K}}\leq v(FD_K)-v(P)=v(FP_K)-v(P)\leq U{\frac{\ln K}{K}}.$$

Furthermore, we will prove $$d((x^*)^K,F(P))\leq \sigma(4nR+3)A\frac{\ln K}{K},$$ where $(x^*)^K\in S(FP_K)$, when $K$ is large enough. If $||S(FP_K)||>2R$, there exists $\{K_i\}$ such that $(x^*)^{K_i}\in S(FP_{K_i})$, $||(x^*)^{K_i}||>2R.$ Taken $(x^{**})^{K_i}\in l(\bar{x},(x^*)^{K_i})$ satisfied $||(x^{**})^{K_i}||=2R,$ where $\bar{x}\in S(P)$, then, $$-L\frac{\ln K_i}{K_i}\leq c^T(x^{**})^{K_i}-v(P)\leq U\frac{\ln K_i}{K_i}.$$ Since $\bar{x}\in F(AFP_{K_i})$ and $(x^{*})^{K_i}\in F(AFP_{K_i})$, then $(x^{**})^{K_i}\in F(AFP_{K_i}).$ By Lemma \ref{lem-0}, $$d((x^{**})^{K_i},F(P))\leq \sigma \max\{\max_{t\in [0,2\pi]} \sum_{j=1}^na_j(t)(x^{**})^{K_i}_j-a_0(t) ,0\}.$$

Noticing that \begin{eqnarray}
              &&\max_{t\in [0,2\pi]} \sum_{j=1}^na_j(t)(x^{**})^{K_i}_j-a_0(t) \nonumber\\
                 =&&  \max_{t\in [0,2\pi]} \sum_{j=1}^n(a_j(t)-\hat{a}_j^{K_i}(t))(x^{**})^{K_i}_j-(a_0(t)-\hat{a}_0^{K_i}(t))+  \sum_{j=1}^n\hat{a}_j^{K_i}(t)(x^{**})^{K_i}_j-\hat{a}_0^{K_i}(t)\nonumber\\
                 \leq &&(2nR+1)A\frac{\ln K_i}{K_i}+2(nR+1)A\frac{\ln K_i}{K_i} \nonumber\\
                 =&&(4nR+3)A\frac{\ln K_i}{K_i}, \nonumber
              \end{eqnarray}
       so $$d((x^{**})^{K_i},F(P))\leq \sigma(4nR+3)A\frac{\ln K_i}{K_i}.$$
When $K_i\rightarrow \infty$, $d((x^{**})^{K_i},F(P)) \rightarrow 0$. For any accumulation point $x^{**}$ of $(x^{**})^{K_i}$, we have  $c^Tx^{**}=v(P),$ $||x^{**}||=2R$ , which is against to  $||S(P)||\leq R$, so $||S(FP_K)||\leq 2R$, when $K$ is large enough.

By Lemma \ref{lem-0}, $\forall (x^*)^K\in S(FP_K)$,
$$d((x^*)^K,F(P))\leq \sigma\max\{\max_{t\in [0,2\pi]}\sum_{j=1}^na_j(t)(x^*)^K_j-a_0(t),0\}.$$
By the same discussion above, we could get
$$d((x^*)^K,F(P))\leq \sigma(4nR+3)A\frac{\ln K}{K},$$
so $(x^*)^K$ is a $\displaystyle\frac{\ln K}{K}-$ optimal solution of $(P)$ when $K$ is large enough.

 (d) We will prove $((x^*)^K,(\bar{Y}^*)^K)\in S(TFP_K)$ with some $(\bar{Y}^*)^K\succeq 0.$ Since the slater condition  holds for Program $(FP_K)$,  then, $v(TFP_K)\geq v(TFD_K)=v(FD_K)=v(FP_K)$, and
 \begin{eqnarray*}
              0=\min & \displaystyle\int_{0}^{2\pi}(\hat{a}^K_0(t)-\sum_{j=1}^n\hat{a}^K_j(t)(x^*)^K_j)\mu(dt) \\
              s.t.&  \mu\in \mathbf{P}([0,2\pi]),
            \end{eqnarray*}
which means $$\left\{
                \begin{array}{c}
                  \displaystyle\sum_{k=-K}^K(a_{0,k}-\sum_{j=1}^na_{j,k}(x^*)^K_j)y_k<0 \\
                   \\  \left(  \begin{array}{cccc}
     y_0 & y_1 & \cdots & y_K \\
     y_{-1} & y_0 & \cdots & y_{K-1} \\
     \vdots & \vdots & \ddots & \vdots \\
     {y}_{-K} & {y}_{-K+1} & \cdots & y_0 \\
   \end{array}\right)\succeq 0
                \end{array}
              \right.
$$
doesn't have solution. Then, there exists $Y^*\succeq 0$ such that $$a_{0,k}-\displaystyle\sum_{j=1}^na_{j,k}(x^*)^K_j=(\mathcal{C}_1^K(Y^*)^K)_k,  k=-K,\cdots,K.$$
So $(\bar{Y}^*)^K=\left(
        \begin{array}{cc}
          (Y^*)^K_r & -(Y^*)_s^K \\
          (Y^*)^K_s & (Y^*)^K_r \\
        \end{array}
      \right)\succeq 0,
$ where $$(Y^*)^K=(Y^*)^K_r+i(Y^*)^K_s, (Y^*)^K_r, (Y^*)^K_s\in \mathbb{R}^{(K+1)\times(K+1)},$$ and  $((x^*)^K,(\bar{Y}^*)^K)\in F(TFP_K),$ which means $v(TFP_K)\leq v(FP_K)$ and  $((x^*)^K,(\bar{Y}^*)^K)\in S(TFP_K).$ \qed

\begin{remark}

 In general, condition (b) in Theorem \ref{thm-0} does not hold for $a_j(t)$. If we replace $a_j(t)$ as
$$\bar{a}_j(t)=\left\{
                 \begin{array}{cc}
                   a_j(2\pi-2t), & 0\leq t\leq \pi \\
                   a_j(2t-2\pi), & \pi<t\leq 2\pi \\
                 \end{array}
               \right.,
$$
then $\bar{a}_j(t)$ is Lipschitz continuous with $\bar{a}_j(0)=\bar{a}_j(2\pi)$ and could be considered as a periodic function with 2$\pi$ period, while the feasible set of  Program $(P)$ is kept. Then, $\bar{a}_j(t)$ must be even function, which means $$\frac{1}{2\pi}\int_{0}^{2\pi}\bar{a}_j(t)e^{ikt}dt$$
must be real, and $r_{j,k}=0$.
 So Program $(TFP_K)$ is equivalent to the following program  denoted as Program ($RTFP_K$):
\begin{subequations}\label{eq-9}
\begin{eqnarray}
  (RTFP_K) ~\min & c^Tx &\nonumber\\
        s.t. & {\mathcal{C}}^K_1Y+(E^K)^Tx+p^K=0 &\nonumber\\
        & Y\succeq 0, \nonumber
\end{eqnarray}
   \end{subequations}
 where $\mathcal{C}^K_1=\{\mathcal{C}^K_{1,j}\}_{j=0}^K.$

\end{remark}
\begin{cor}\label{cor-1}
If the conditions (a) and (b) in Theorem \ref{thm-0} are replaced by
 \begin{itemize}
   \item[(a')] $a_j(t)$ has $q-$th  derivative $a_j^{(q)}(t)$ and  $a_j^{(q)}(t)$ are Lipschitz continuous functions on $[0, 2\pi]$;
   \item[(b')]  $a_j^{(q)}(0)=a_j^{(q)}(2\pi)$, $j=0,\cdots,n$;
 \end{itemize}
 respectively, $x^{*}$ is a $\displaystyle\frac{\ln K}{K^{q+1}}-$optimal solution of Program $(P)$,   when $K$ is large enough, where $((x^*)^K,(Y^*)^K)\in S(TFP_K)$ with some $(Y^*)^K\succeq 0.$
\end{cor}
{\noindent{\rm{ \bf Proof~}} } The proof procedure is similar to that of Theorem  \ref{thm-0}, except that  we use  Corollary III in page 22 of \cite{book2} instead of Corollary I. $\hfill \square$

	As is evident from the foregoing discussion, in order to solve Program $(TFP_K)$ or $(RTFP_K)$, both of which almost have the same parameters $a_{j,k}$, except the case $j=0,k=0$, the critical step in this algorithmic process is to efficiently compute the parameters $a_{j,k}$.
Noting that
$$a_{j,k} =\frac{1}{2\pi}\int_{0}^{2\pi}a_j(t)e^{ikt}dt,\ k=-K,\ldots,K, $$ which is the $k$th coefficient of the Fourier transform of $a_j(t)$, we are only required to calculate $a_{j,k}$, for $k\geq 0$. In a general setting, we cannot accurately obtain $a_{j,k}$, and hence, we resort to using the {\it discrete Fourier transformations} ({\bf DFT}) to compute the necessary coefficients $a_{j,k}$  in a relatively quick fashion. The {\it fast Fourier transform} ({\bf FFT}) is a well-known algorithm for computing discrete Fourier transforms, or their inverse \cite{FFT}. A direct application of the DFT algorithm towards determining the coefficients $a_{j,k}$  results in a complexity of the order $O(K^2)$. But by exploiting the symmetry and periodicity of $e^{ikt}$, it is feasible to reduce this complexity to $O(K\ln K)$, where $K$ is the truncation factor, which is a measure of the accuracy of the approximation and the data size involved in the problem.

Here, we present our algorithm.
 \begin{alg}\label{alg-1}
    \item[Step 0:] Given $\varepsilon>0$, $N>0$ large enough.  Set
      $$\bar{a}_j(t)=\left\{
                       \begin{array}{cc}
                        a_j(2\pi-2t), & 0\leq t\leq \pi \\
                         a_j(2t-2\pi), & \pi<t\leq 2\pi \\
                       \end{array}
                     \right., \ j=0,\cdots,n;
      $$
    \item[Step 1:] Calculate $\bar{a}_{j,k}, \ k=0,\cdots,N, j=0,\cdots,n$ (by FFT), and get the real part $r_{j,k}$  and the imaginary part $s_{j,k}$ of them;
    \item[Step 2:] Solving Program $(TFP_K)$ or $(RTFP_K)$  to get the approximate optimal solution of Program $(P)$.
 \end{alg}

 \begin{remark} There are two variants of Algorithm \ref{alg-1} that have been used in our computational implementations, one that uses the FFT method and the other that does not, depending on the problem structure. Specifically, if the frequency of $a_j(t)$ is not too large, we can gainfully apply the FFT variant to generate good approximations of $a_{j,k}$, but the resulting $s_{j,k}$, i.e., the coefficients of the complex terms in the problem, may not be equal to zero. In this case, we would need to contend with solving the larger Program $(TFP_K)$ as opposed to merely solving Program $(RTFP_K)$; Alternately, we can generate the values of $a_{j,k}$ by using other (perhaps more computationally intensive) schemes, but gain the advantage of solving the resulting relatively smaller Program $(RTFP_K)$.
\end{remark}

\section{Numerical Results}\label{sec5}
In this section, we demonstrate the efficiency of Algorithm \ref{alg-1} by solving five illustrative numerical examples taken from the literature \cite{sip3,lsip3}. As the default index sets of these examples is not in the range $[0,2\pi]$, we make some minor modifications to scale these programs into the required format. All of our computations are conducted on a Windows machine, equipped with a dual core 2.69GHz processor and 8GB RAM, using MATLAB R2011b as the computational engine. For solving the underlying semidefinite programs, we use SDPT4 and YALMIP software packages \cite{yalmip,SDPT3}. These five examples are given below.

\begin{exam}\label{ex-1}
\begin{eqnarray}\label{eq-11}
  \min & \displaystyle \sum_{j=1}^n \frac{1}{j}x_j \nonumber\\
   s.t.& \displaystyle \sum_{j=1}^n -(\frac{t}{2\pi})^{j-1}x_j\leq -\tan(\frac{t}{2\pi}),& \forall t\in [0, 2\pi].
\end{eqnarray}
\end{exam}

\begin{exam}\label{ex-2}
\begin{eqnarray}\label{eq-12}
  \min & \displaystyle \sum_{j=1}^nx_j \nonumber\\
   s.t.& \displaystyle \sum_{j=1}^n -(\frac{t}{2\pi}+1)^{j-1}x_j\leq -\frac{2\pi}{\sqrt{4\pi^2+t^2}},& \forall t\in [0, 2\pi].
\end{eqnarray}
\end{exam}

\begin{exam}\label{ex-3}
\begin{eqnarray}\label{eq-13}
  \min & \displaystyle \sum_{j=1}^8 \frac{1}{j}x_j \nonumber\\
   s.t.& \displaystyle \sum_{j=1}^8 -(\frac{t}{2\pi})^{j-1}x_j\leq -\displaystyle\frac{2\pi}{4\pi-t},& \forall t\in [0, 2\pi].
\end{eqnarray}
\end{exam}
\begin{exam}\label{ex-4}
\begin{eqnarray}\label{eq-14}
  \min & \displaystyle \sum_{j=1}^9 \frac{1}{j}x_j \nonumber\\
   s.t.& \displaystyle \sum_{j=1}^9 -(\frac{t}{2\pi})^{j-1}x_j\leq -\frac{4\pi^2}{{4\pi^2+t^2}},& \forall t\in [0, 2\pi].
\end{eqnarray}
\end{exam}
\begin{exam}\label{ex-5}
\begin{eqnarray}\label{eq-15}
  \min & \displaystyle \sum_{j=1}^{10} -(0.95)^{2j-1}x_j \nonumber\\
   s.t.& \displaystyle \sum_{j=1}^{10} -2\cos(\frac{(2j-1)t}{2})x_j\leq 1,& \forall t\in [0, 2\pi].
\end{eqnarray}
\end{exam}

We begin our computations by solving Example \ref{ex-1} and \ref{ex-2}, corresponding to different values of $n \in \{5, 6, 7, 8\}$, where $n$ denotes the dimension of $x$; $K$ is the truncation factor in Program $(TFP_K)$, which is also the number of Fourier terms presented in the program; $\nu_{Alg 1}$ denotes the  optimal value obtained by Algorithm \ref{alg-1}; and the relative error is the absolute difference between the actual optimal value and the objective function value obtained by our algorithm.
From \cite{sip5,lsip3}, the optimal values of Example \ref{ex-1}, Example \ref{ex-2} for different dimensions  are shown in Table \ref{tab-1}. And by \cite{sip3}, the optimal values of Example \ref{ex-3}, Example \ref{ex-4} and  Example \reff{ex-5} are $0.69314815$, $0.78549953$, and $-0.48354840$, respectively.

\begin{table}[h]
 \caption{Optimal values of Example \ref{ex-1} and Example \ref{ex-2}, corresponding to different dimensions}\label{tab-1}
\begin{center}\normalsize
\begin{tabular}{|c|c|c|c|c|}
  \hline
    n &    5 &  6 & 7 & 8\\
     \hline
    Example \ref{ex-1}&  0.61740424 & 0.61608515&0.61572945&0.61565322\\
    \hline
    Example \ref{ex-2}& 1 &1 &1&1\\
    \hline
\end{tabular}

 \end{center}

\end{table}

%

We give the results of our algorithm in following tables. We use the FFT type of Algorithm \ref{alg-1} for Example \ref{ex-1}, Example \ref{ex-2},  Example \ref{ex-3},  and Example \ref{ex-4}, since the frequencies of all the functions $a_j(t)$ are not too large. For Example \ref{ex-5}, we use Algorithm \ref{alg-1} without FFT. For comparing, we use the cutting plant method \cite{sip5} to compute Example \ref{ex-1} and \ref{ex-2}. For these two examples, we test 10 times and get the average results showed in Table \ref{tab-3.5} by using the cutting plant method.

\begin{table}[h]
     \caption{Numerical results by using Algorithm \ref{alg-1} on Example \ref{ex-1}} \label{tab-2}
\begin{center}\normalsize
  \begin{tabular}{|c|c|c|c|c|c|c|c|c|c|}
    \hline
     &   \multicolumn{3}{c|}{$\nu_{\text{Alg \ref{alg-1}}}$}  & \multicolumn{3}{c|}{\text{CPU Time (seconds)}}  &\multicolumn{3}{c|}{\text{Relative Error}}  \\
      \hline
 n $\diagdown$ K & 8 & 16 & 32 & 8 & 16 & 32 & 8 & 16 & 32 \\
     \hline
    5 & 0.6151 & 0.6164 &  0.6170& 0.323 & 0.552 & 1.994 & $2.3\times10^{-3}$ &$ 1.0\times10^{-3}$& $ 4.0\times10^{-4}$ \\
     \hline
    6 & 0.6152 &  0.6158 &  0.6159 & 0.305 &  0.620 & 2.224 & $9.0\times10^{-4}$ &$3.0\times10^{-4}$ & $1.0\times10^{-4}$ \\
     \hline
    7 &  0.6154 &  0.6156 & 0.6157 & 0.326 & 0.694 & 2.44 & $3.0\times10^{-4}$ & $1.0\times10^{-4}$ & $3.0\times10^{-5}$\\
     \hline
    8 &  0.6158 & 0.6156 &  0.6156 &  0.401 & 0.577 &  2.515 & $1.9\times10^{-3}$ & $2.0\times10^{-4}$ & $1.0\times10^{-4}$ \\
    \hline
  \end{tabular}

   \end{center}

    \end{table}
    \begin{table}[h]
     \caption{Numerical results by using Algorithm \ref{alg-1} on Example \ref{ex-2}} \label{tab-3}
    \begin{center}\normalsize
  \begin{tabular}{|c|c|c|c|c|c|c|c|c|c|}
    \hline
     &   \multicolumn{3}{c|}{$\nu_{\text{Alg \ref{alg-1}}}$}  & \multicolumn{3}{c|}{\text{CPU Time (seconds)}}  &\multicolumn{3}{c|}{\text{Relative Error}}  \\
      \hline
  n $\diagdown$ K & 8 & 16 & 32 & 8 & 16 & 32 & 8 & 16 & 32 \\
     \hline
    5 & 1.00 & 1.00 &  1.00&  0.388 &  0.774 & 3.397 & $2.7\times10^{-6}$ &$ 1.6 \times10^{-6}$& $ 3.0\times10^{-7}$ \\
     \hline
    6 & 1.00 &  1.00 &  1.00 &  0.502 &  0.843 &  3.305 & $1.3\times10^{-6}$ &$8.3\times10^{-9}$ & $1.8\times 10^{-8}$ \\
     \hline
    7 &  1.00 &  1.00&1.00 & 0.397 & 0.775 &  3.151 & $2.0 \times10^{-7}$ & $2.0 \times10^{-7}$ & $4.9\times 10^{-8}$\\
     \hline
    8 &  1.00 & 1.00 &  1.00 &   0.427 & 0.789 &   2.864 & $3.0\times10^{-7}$ & $1.0\times10^{-7}$ & $1.2\times 10^{-8}$ \\
    \hline
  \end{tabular}

   \end{center}
\end{table}

    \begin{table}[h]
    \begin{center}\normalsize
     \caption{Numerical results by using the cutting plant method  on Example \ref{ex-1} and Example \ref{ex-2}} \label{tab-3.5}
  \begin{tabular}{|c|c|c|c|}
    \hline
      n & Example \ref{ex-1} & \text{CPU Time (seconds)} & \text{Relative Error}\\
      \hline
     5 &    0.6165                  &  0.437                         & $9.0\times 10^{-4}$    \\
     \hline
     6 &     0.6154                 &     0.520                      &  $6.0\times 10^{-4}$   \\
     \hline
     7 &    0.6155                  & 0.598                        &   $5.5\times 10^{-4}$   \\
     \hline
     8 &     0.6154                 &       0.588                    &  $7.0\times 10^{-4}$ \\
    \hline
    \end{tabular}
     \begin{tabular}{|c|c|c|c|}
     \hline
   n &Example \ref{ex-2} & \text{CPU Time (seconds)}& \text{Relative Error} \\
     \hline
     5 &         0.999             &          0.245                 &  $3.0\times 10^{-5}$    \\
     \hline
     6 &         1.000             &        0.334                   &  $8.0\times 10^{-6}$   \\
     \hline
     7 &         1.000             &      0.350                     &  $7.2\times 10^{-6}$    \\
     \hline
     8 &         1.000             &       0.421                    &  $4.3\times 10^{-6}$  \\
    \hline
  \end{tabular}

   \end{center}
\end{table}

 From the results recorded in Tables \ref{tab-2}, \ref{tab-3}, and \ref{tab-3.5}, we find that our algorithm has  better precision than the cutting plant method, although our times need a little more. But from the numerical tests, we found the results of the cutting plants are strongly depended on the start scattered point, which means that the optimal points are quite different when using the different start scattered  point, while ours are more stable. Furthermore, from Tables \ref{tab-2} and \ref{tab-3}, we infer that using a truncation factor $K \approx 2n$ yields an acceptable approximation to the original program (with objective function value within the prescribed tolerance limit), while yet expending only a reasonable computational effort.

The results for Example \ref{ex-3}, Example \ref{ex-4} and Example \ref{ex-5} are shown in Table \ref{tab-4} respectively, where $K=20$.
\begin{table}[h]
    \caption{Numerical results by using Algorithm \ref{alg-1} on Example \ref{ex-3}, Example \ref{ex-4}, and Example \ref{ex-5}}\label{tab-4}
\begin{center}\normalsize
    \begin{tabular}{|c|c|c|c|}
      \hline
     Example   &$\nu_{\text{Alg \ref{alg-1}}}$ & \text{CPU Time (seconds)} & \text{Relative Error}\\
      \hline
   \ref{ex-3} & 0.6931 & 0.942&  $7.7\times 10^{-7}$\\
      \hline
   \ref{ex-4} &  0.7854 &  0.780 & $2.4\times 10^{-4}$ \\
      \hline
   \ref{ex-5}& -0.4835 &  0.411  &$1.2\times 10^{-9}$\\
      \hline
    \end{tabular}

\end{center}

\end{table}

 We could see that our algorithm is efficient and effective.
\section{Conclusions}\label{sec6}

In this paper, we present a framework to determine  solutions for a class of linear semi-infinite programs  based on the theory of trigonometric moments. Beginning with the dual program to the linear semi-infinite program, we first construct an approximation to the dual program using Fourier transformations, and subsequently transform the truncated Fourier dual program into an approximate semidefinite program. Recognizing that the derived semidefinite program is not a standard optimization problem due to the presence of complex variables, we cleverly exploit the properties of the underlying matrices to derive an equivalent real-valued semidefinite program. Moreover, we also discuss various relationships between the original semi-infinite program and the derived semidefinite program, and prove that the optimal solution of the truncated semidefinite program does indeed be the approximate optimal solution of the original semi-infinite program. An algorithm based on Fast Fourier Transforms is presented, and our illustrative computational results serve to showcase the computational efficiency of the algorithm.

\bigskip\medskip
\noindent{\bf  Acknowledgements:} \small{Y. Xu was supported in part by National Natural Science Foundation of China No. 11501100,  11571178, 11671082 and 11871149. J. Desai was supported in part by the Ministry of Education (Singapore) Academic Research Fund Tier 1 Grant No. M4011083. X. Yan was supported in part by the STIP of Higher Education Institutions in Shanxi No. 201802103, and National Natural Science Foundation of China No. 11901424.
}


\begin{thebibliography}{abc99xyz}\normalsize

\bibitem{trm1} Akhiezer, N.I.: The classical moment problem. Hafner, New York (1965)

\bibitem{trm2}  Akhiezer, N.I., Krein, M.G.: Some questions in the theory of moments. American
Mathematical Society Translations Vol. 2 (1962)

\bibitem{sip3}  Betr\`{o}, B.: An accelerated central cutting plane algorithm for linear semi-infinite programming.
Math. Program., 101, 479-495 (2004)


\bibitem{sip5}  Goberna, M.A., L\'{o}pez, M.A.: Linear semi-infinite optimization. Chichester: John Wiley \& Sons
(1998)

\bibitem{sipb1} Hettich, R., Kortanek, K.O.:  Semi-infinite programming: theory, methods, and applications. SIAM
Rev., 35:380-429 (1993)


\bibitem{lmb} Lasserre, J.B.: Moments, positive polynomials and their applications. Imperial
College Press Optimization Series, vol. 1, Imperial College Press, London (2010)

\bibitem{lmp} Lasserre, J.B.: Global optimization with polynomials and the problem
of moments. SIAM Journal on Optimization, 11(3), 796-817 (2001)

\bibitem{book2} Jackson, D: The theory of approximation. American Mathematical Society, New York (1930)

\bibitem{sipb3} Ling, C., Ni, Q., Qi, L., Wu, S.Y.: A new smoothing Newton-type algorithm for semi-infinite
programming. J. Glob. Optim., 47, 133-159 (2010)

\bibitem{yalmip} L\"{o}fberg, J.: YALMIP : A Toolbox for Modeling and Optimization in MATLAB. In Proceedings of the CACSD Conference, Taipei, Taiwan (2004)

\bibitem{LuoTseng} Luo, Z-Q, Tseng, P.: Pertubation analysis of a condition number for linear systems, SIAM J. Matrix Anal. Appl.
Vol. 15, No. 2, 636-660,  (1994)

\bibitem{sipb2} Reemtsen, R., R\"{u}ckmann, J-J, ,editors: Semi-infinite programming. Boston: Kluwer Academic
(1998)

\bibitem{sip4}  Stein, O., Still, G.: Solving semi-infinite optimization problems with interior point techniques.
SIAM J. Control Optim., 42, 769-788 (2003)

\bibitem{SDPT3}  Tutuncu, R.H., Toh, K.C., Todd, M.J.: Solving
semidefinite-quadratic-linear programs using SDPT3.  Math.
Prog., 95, 189-217 (2003)

\bibitem{FFT} van Loan, C.: Computational Frameworks for the Fast Fourier Transform.
Philadelphia: SIAM (1992)

\bibitem{book-sdp1} Wolkowicz, H., Saigal, R., Vandenberghe, L.: Handbook of
Semidefinite Programming. Kluwer (2000)


\bibitem{lsip3} Xu, Y., Sun, W.Y., Qi, L.: On solving a class of linear semi-infinite
programming by SDP method. Optimization,  64(3), 603-616 (2015)


\end{thebibliography}
\end{document}